\documentclass[12pt]{article}
\usepackage{fullpage,graphicx,graphics,amsmath,amsfonts,verbatim,mathtools,subcaption}
\usepackage{url,bbm}
\usepackage{amssymb,enumitem}
\usepackage{booktabs}

\newcommand{\downto}{\downarrow}

\newcommand{\BEAS}{\begin{eqnarray}}
\newcommand{\EEAS}{\end{eqnarray}}
\newcommand{\BEA}{\begin{eqnarray}}
\newcommand{\EEA}{\end{eqnarray}}
\newcommand{\BEQ}{\begin{equation}}
\newcommand{\EEQ}{\end{equation}}
\newcommand{\BIT}{\begin{itemize}}
\newcommand{\EIT}{\end{itemize}}
\newcommand{\BNUM}{\begin{enumerate}}
\newcommand{\ENUM}{\end{enumerate}}

\newcommand{\BA}{\begin{array}}
\newcommand{\EA}{\end{array}}

\newcommand{\eg}{{\it e.g.}}
\newcommand{\ie}{{\it i.e.}}

\newcommand{\ones}{\mathbf 1}

\newcommand{\reals}{{\mbox{\bf R}}}

\newcommand{\argmin}{\mathop{\rm argmin}}

\newcommand{\dom}{\mathop{\bf dom}}

\newcounter{algorithmctr}[section]
\renewcommand{\thealgorithmctr}{\thesection.\arabic{algorithmctr}}
\newenvironment{algdesc}%
   {\refstepcounter{algorithmctr}\begin{list}{}{%
       \setlength{\rightmargin}{0\linewidth}%
       \setlength{\leftmargin}{.05\linewidth}}%
       \rmfamily\small
       \item[]{\setlength{\parskip}{0ex}\hrulefill\par%
        \nopagebreak{\bfseries\textsf{Algorithm \thealgorithmctr~}}}}%
   {{\setlength{\parskip}{-1ex}\nopagebreak\par\hrulefill} \end{list}}

\title{Minimizing a Sum of Clipped Convex Functions}
\author{Shane Barratt \and Guillermo Angeris \and Stephen Boyd}

\begin{document}
\maketitle

\begin{abstract}
We consider the problem of minimizing a sum of clipped convex functions;
applications include clipped empirical risk minimization
and clipped control.
While the problem of minimizing the sum of clipped convex functions is NP-hard, we present some
heuristics for approximately solving instances of these problems.
These heuristics can be used to find good, if not global, solutions and appear to work well in practice.
We also describe an alternative formulation, based on the perspective transformation,
which makes the problem amenable to mixed-integer convex programming
and yields computationally tractable lower bounds.
We illustrate one of our heuristic methods by applying it to various examples and use the perspective transformation to certify that the solutions are relatively close to the global optimum.
This paper is accompanied by an open-source implementation.
\end{abstract}

\section{Introduction}
Suppose $f: \reals^n \to \reals$ is a convex function, and 
$\alpha \in \reals$.   We refer to the function $\min\{f(x),\alpha\}$
as a \emph{clipped convex function}.
In this paper we consider the problem of minimizing a sum of 
clipped convex functions,
\BEQ
\begin{array}{ll}
\mbox{minimize} & f_0(x) + \sum_{i=1}^m \min\{f_i(x),\alpha_i\},
\end{array}
\label{eq:main_formulation}
\EEQ
with variable $x\in\reals^n$, where
$f_0:\reals^n \to \reals \cup \{+\infty\}$
and $f_i:\reals^n \to \reals$ for $i=1,\ldots,m$ are closed proper convex functions,
and $\alpha_i\in\reals$ for $i=1,\ldots,m$.
We use infinite values of $f_0$ to encode constraints on $x$, \ie, to
constrain $x \in \mathcal X$ for a closed convex set $\mathcal X$
we let $f_0(x)=+\infty$ for all $x \not \in \mathcal X$.
When $f_i(x) > \alpha_i$, the value of the $i$th term in the sum
is \emph{clipped} to $\alpha_i$, which limits how large each term 
in the objective can be.
Many practical problems can be formulated as instances 
of~\eqref{eq:main_formulation};
we describe a few in~\S\ref{sec:applications}.

\paragraph{NP-hardness.}
In general, problem~\eqref{eq:main_formulation} is nonconvex and as a result
can be very difficult to solve.
Indeed,~\eqref{eq:main_formulation} is NP-hard.
We show this by giving a reduction
of the subset sum problem to an instance of~\eqref{eq:main_formulation}.

The subset sum problem involves determining whether or not there exists a subset
of a given set of integers $a_1,\ldots,a_n$ that sum to zero.
The optimal value of the problem
\begin{equation*}
\begin{array}{ll}
\mbox{minimize} & (a^Tx)^2 - n/4 + \sum_{i=1}^n \min\{x_i^2, 1/4\} + \min\{(x_i-1)^2, 1/4\} \\
\mbox{subject to} & \ones^T x \geq 1,
\end{array}
\end{equation*}
which has the form~\eqref{eq:main_formulation},
is zero if and only if $x_i \in \{0, 1\}$, at least one of $x_i=1$,
and $a^T x = 0$; in other words, the set $\{a_i \mid x_i = 1\}$ sums to zero.
Since the subset sum problem can be reduced to
an instance of~\eqref{eq:main_formulation}, we conclude that in general
our problem is at least as hard as difficult problems
like the subset sum problem.

\paragraph{Global solution.}
There is a simple (exhaustive) method to solve~\eqref{eq:main_formulation} globally:
for each subset $\Omega$ of $\{1,\ldots,m\}$, we solve the convex problem
\begin{equation}\label{eq:convex-subset}
\begin{array}{ll}
\mbox{minimize} & f_0(x) + \sum_{i \in \Omega} f_i(x) + \sum_{i \not \in \Omega} \alpha _i \\
\mbox{subject to} & f_i(x) \leq \alpha_i, \quad i \in \Omega,
\end{array}
\end{equation}
with variable $x \in \reals^n$.
The solution to~(\ref{eq:convex-subset}) with the lowest optimal
value is the solution to~\eqref{eq:main_formulation}.
This general method is not practical unless $m$ is quite small,
since it requires the solution of $2^m$ convex optimization problems.

In some specific instances of problem~\eqref{eq:main_formulation}, 
we can cut down the search space if we know that a specific 
choice of $\Omega \subseteq \{1, \dots, m\}$ implies
\[
\{x \mid f_i(x) \leq \alpha_i,\, i \in \Omega\} = \emptyset,
\]
which means that the optimal value of~(\ref{eq:convex-subset}) is
$+\infty$.
In this case, we do not have to 
solve problem~\eqref{eq:convex-subset} for this choice of $\Omega$, as we know it will be infeasible.
One simple example where this happens is when the $\alpha_i$-sublevel sets of $f_i$ are pairwise disjoint, which implies that
we only have to solve $m$ convex problems (as opposed to $2^m$) to find the global solution.
This idea is used in~\cite{minimizing2019liu} to guide their proposed search algorithm.

\paragraph{Related work.}
The general problem of minimizing a sum of clipped convex functions
was recently considered in~\cite{minimizing2019liu}.
In their paper, they also show that the problem is NP-hard via a reduction to 3-SAT and
give a global solution method in a few special cases whenever $n$ is small.
They also provide a heuristic method based on cyclic coordinate descent,
leveraging the fact that one-dimensional problems are easy to solve.

The idea of using clipped convex functions has appeared
in multiple application areas, the most prominent being statistics.
For example,
the sum of clipped absolute values (often referred to as the \emph{capped} $\ell_1$-norm)
has been used as a sparsity-inducing regularizer
\cite{zhang2009multi, zhang2010analysis, ong2013learning}.
In particular,~\cite{zhang2009multi, ong2013learning}
make use of the fact that problem~\eqref{eq:main_formulation}
can be written as a difference-of-convex (DC) problem and can be
approximately minimized via the convex-concave procedure \cite{lipp2016variations}
(see Appendix~\ref{sec:convex_concave}).
The clipped square function (also known as the \emph{skipped-mean} loss)
was also used in~\cite{torr1998robust} to estimate
view relations, and in~\cite{portilla2015efficient} to
perform robust image restoration.
Similar approaches have been taken for clipped loss functions,
where they have been used for robust feature selection~\cite{lan2016robust},
regression~\cite{yang2010relaxed,she2011outlier},
classification~\cite{suzumura2014outlier,safari2014insensitive,xu2016robust},
and robust principal component analysis~\cite{sun2013robust}.

\paragraph{Summary.}
We begin by presenting some applications of minimizing
a sum of clipped convex functions in~\S\ref{sec:applications}
to empirical risk minimization and control.
We then provide some simple heuristics for approximately
solving~\eqref{eq:main_formulation} in~\S\ref{sec:methods},
which we have found to work well in practice.
In~\S\ref{sec:perspective_formulation}, we describe a method for
converting~\eqref{eq:main_formulation} into a mixed-integer convex program,
which is amenable to solvers for mixed-integer convex programs.
Finally, we describe an open-source Python implementation of the ideas
described in this paper
in~\S\ref{sec:implementation} and apply our implementation to a
few illustrative examples in~\S\ref{sec:examples}.

\section{Applications}
\label{sec:applications}
In this section we describe some possible
applications of minimizing a sum of clipped convex functions.

\subsection{Clipped empirical risk minimization}
\label{sec:clipped_erm}
Suppose we have data
\[
x_1,\ldots,x_N \in \reals^n, \quad y_1,\ldots,y_N \in \mathcal Y.
\]
Here $x_i$ is the $i$th feature vector, $y_i$ is its corresponding output (or label),
and $\mathcal Y$ is the output space.

We find parameters $\theta\in\reals^n$ of a linear model given the data
by solving the \emph{empirical risk minimization} (ERM) problem
\begin{equation}
\begin{array}{ll}
\label{eq:erm}
\mbox{minimize} & \frac{1}{N}\sum_{i=1}^N l(x_i^T\theta,y_i) + r(\theta),
\end{array}
\end{equation}
with variable $\theta$, where $l:\reals \times \mathcal Y \to \reals$ is the loss function,
and $r:\reals^n \to \reals$ is the regularization function.
Here the objective is composed of two parts:
the loss function, which measures the accuracy of the predictions,
and the regularization function, which measures the complexity of $\theta$.
We assume that $l$ is convex in its first argument and that $r$ is convex,
so the problem~(\ref{eq:erm}) is a convex optimization problem.

For a given $x\in\reals^n$, our prediction of $y$ is
\[
\hat y = \underset{y \in \mathcal Y}{\argmin} \;  l(x^T\theta^\star, y),
\]
where $\theta^\star$ is optimal for~(\ref{eq:erm}).
For example, in linear regression, $\mathcal Y = \reals$, $l(z, w)=(z - w)^2$,
and $\hat y = x^T \theta^\star$;
in logistic regression, $\mathcal Y = \{-1,1\}$, $l(z, w)=\log(1+e^{-wz})$,
and $\hat y = \mathbf{sign}(x^T\theta^\star)$, 
where $\mathbf{sign}(z)$ is equal
to $1$ if $z \geq 0$ and $-1$ otherwise.

While ERM often works well in practice, 
it can perform poorly when there are outliers in the data.
One way of fixing this is to clip the loss for each data point
to a value $\alpha\in\reals$, leading to the \emph{clipped ERM} problem,
\begin{equation}
\label{eq:clipped-erm}
\begin{array}{ll}
\mbox{minimize} & \frac{1}{N}\sum_{i=1}^N \min\{l(x_i^T\theta,y_i),\alpha\} + r(\theta).
\end{array}
\end{equation}
After solving (or approximately solving) the clipped problem,
we can label data points $(x_i,y_i)$ where $l(x_i^T\theta^\star,y_i) \geq \alpha$
as outliers.
The clipped ERM problem is an instance
of what is referred to in statistics as a \emph{redescending M-estimator}
\cite[\S 4.8]{huber2009robust},
since the derivative of the clipped loss goes to $0$ as the magnitude of its input goes to infinity.
In this terminology, the clip value $\alpha$ is referred to
as the \emph{minimum rejection point}.

In~\S\ref{sec:erm-example}, we show an example where the normal empirical risk minimization problem fails, while its clipped variant has 
good performance.

\subsection{Clipped control}
Suppose we have a linear system with dynamics given by
\[
x_{t+1} = Ax_t + Bu_t, \quad t=0,\ldots,T-1,
\]
where $x_t \in \reals^n$ is the state of the system
and $u_t\in\reals^p$ denotes the input to the system, at time period $t$.
The dynamics matrix $A \in \reals^{n \times n}$ and
the input matrix $B \in \reals^{n \times m}$ are given.

We are given stage cost functions $g_t:\reals^n \times \reals^p \to \reals$,
and an initial state $x^\mathrm{init}\in\reals^n$.
The standard optimal control problem is
\[
\begin{array}{ll}
\mbox{minimize} & \sum_{t=0}^{T} g_t(x_t, u_t)\\
\mbox{subject to} & x_{t+1} = A_t x_t + B_t u_t, \quad t=0,\ldots,T-1, \\
& x_t\in\mathcal X_t, \quad u_t\in\mathcal U_t, \quad t=0,\ldots,T,\\
& x_0 = x^\mathrm{init},
\end{array}
\]
where, at time $t$, $\mathcal X_t\subseteq\reals^n$ is the convex set of allowable states
and $\mathcal U_t\subseteq\reals^m$ is the convex set of allowable inputs.
The variables in this problem are the states and inputs,
$x_t$ and $u_t$.
If the stage cost function $g_t$ are convex, the optimal control
problem is a convex optimization problem.

We define a \emph{clipped optimal control} problem as
an optimal control problem in which
the stage costs can be expressed as sums of clipped convex functions,
\ie,
\[
g_t(x,u) = g_t^0(x,u) + \sum_{i=1}^K \min\{g_t^{i}(x,u),\alpha_t^i\},
\]
where, for all $t$ and $i=1,\ldots,K$, the functions $g_t^{i}:\reals^n\times\reals^m\to\reals$ are convex
and $\alpha_t^i\in\reals$.
This gives another instance of our general 
problem~(\ref{eq:main_formulation}).

A simple but practical example of a clipped control problem is
described in~\S\ref{sec:control-example}.  The problem is to 
design a lane change trajectory for a vehicle; the stage
cost is small when the vehicle is centered in either lane,
which we express as a sum of two clipped convex functions.

\section{Heuristic methods}
\label{sec:methods}
There are many methods for
approximately solving~\eqref{eq:main_formulation}.
In this section we describe a few heuristic methods
that we have observed to work well in practice.

\paragraph{Bi-convex formulation.}
Throughout this section, we will make use of a simple reformulation
of~\eqref{eq:main_formulation} as the bi-convex problem
\BEQ
\begin{array}{ll}
\mbox{minimize} & L(x,\lambda)
= f_0(x) + \sum_{i=1}^m \lambda_i f_i(x) + (1-\lambda_i)\alpha_i\\
\mbox{subject to} & 0 \le \lambda \le \ones,
\label{eq:nonlinear}
\end{array}
\EEQ
with variables $\lambda \in \reals^m$ and $x \in \reals^n$.
(We note that this reformulation was also pointed out in~\cite[\S3]{yang2010relaxed}.)
The equivalence follows immediately
from the fact that
\[
\min\{a, b\} = \min_{0 \le \lambda \le 1} \left(\lambda a + (1- \lambda) b\right).
\]

\paragraph{Nonlinear programming.}
When $f_i$ are all smooth functions
and $\dom f_0$ is representable as the sublevel set of a smooth function,
it is possible to use general nonlinear solvers to (approximately)
solve~\eqref{eq:nonlinear}.

\paragraph{Alternating minimization.}
Another possibility is to perform alternating minimization on~\eqref{eq:nonlinear},
since each respective minimization is a convex optimization problem.
In alternating minimization, at iteration $k$,
we solve~\eqref{eq:nonlinear} while fixing $\lambda=\lambda^{k-1}$,
resulting in $x^k$.
We then solve~\eqref{eq:nonlinear} while fixing $x=x^k$,
resulting in $\lambda^k$.
It can be shown that
\BEQ
\label{eq:lam_update}
(\lambda^k)_i = \begin{cases}
1 & f_i(x^k) \leq \alpha_i \\
0 & \text{otherwise},
\end{cases}
\EEQ
is a solution for minimization over $\lambda$ with fixed $x = x^k$.

\paragraph{Inexact alternating minimization.}
Although alternating minimization often works well, we have found that
inexact minimization over $\lambda$ works better in practice.
Instead of fully minimizing over $\lambda$, we instead
compute the gradient of the objective with respect to $\lambda$,
\[
g_i = (\nabla_\lambda L(x^k, \lambda))_i = f_i(x^k) - \alpha_i.
\]
We then perform a signed projected gradient
step on $\lambda$ with a fixed step size $\beta > 0$ (we have found $\beta=0.1$ works well in practice, though a range of values all appear to work equally as well). This results in the update
\[
\lambda^{k} = \Pi_{[0,1]^m}(\lambda^k - \beta \mathbf{sign}(g)),
\]
where $\mathbf{sign}$ is applied elementwise to $g$,
and $\Pi_{[0,1]^m}$ denotes the projection onto the unit box, given by
\[
(\Pi_{[0,1]^m}(z))_i = \begin{cases}
1 & z_i \geq 1, \\
z_i & 0 < z_i < 1, \\
0 & \text{otherwise}.
\end{cases}
\]
The final algorithm is described below.
\begin{algdesc}
\label{alg:twostep}
\emph{Inexact alternating minimization.}
\begin{tabbing}
    {\bf given} initial $\lambda^0=(1/2)\ones$, step size $\beta=0.1$, and tolerance $\epsilon > 0$.\\
    {\bf for} $k=1,\ldots,n_\mathrm{iter}$\\
        \qquad \=\ 1.\ \emph{Minimize over $x$.} Set $x^k$ to the solution of the problem\\
            $
            \hspace*{3.5cm} \begin{array}{ll}
            \mbox{minimize} & f_0(x) + \sum_{i=1}^m \lambda^{k-1}_i f_i(x) + (1 - \lambda^{k-1}_i)\alpha_i.
            \end{array}
            $
            \\
        \qquad \=\ 2.\ \emph{Compute the gradient.} Set $g_i=f_i(x^{k}) - \alpha_i$. \\
        \qquad \=\ 2.\ \emph{Update $\lambda$.} Set
            $\lambda^k = \Pi_{[0,1]^m}(\lambda^{k-1} - \beta \mathbf{sign}(g))$. \\
        \qquad \=\ 3.\ \emph{Check stopping criterion.} Terminate if $\|\lambda^k - \lambda^{k-1}\|_1 \le \epsilon$.\\
    {\bf end for}
\end{tabbing}
\end{algdesc}
Algorithm~\ref{alg:twostep} is a descent algorithm in the sense that the
objective function of~\eqref{eq:nonlinear}
decreases after every iteration.
It is also guaranteed to terminate in a finite amount of time,
since there is a finite number of possible values of $\lambda$.
We also note that alternating minimization can be thought of
as a special case of algorithm~\ref{alg:twostep} where $\beta\geq1$.
In practice, we have found that algorithm~\ref{alg:twostep}
often finds the global optimum in simple problems and appears
to work well on more complicated cases.
We use algorithm~\ref{alg:twostep} in our
generic \verb|cvxpy| implementation (see \S\ref{sec:implementation}).

\section{Perspective formulation}
\label{sec:perspective_formulation}

In this section we describe the perspective formulation of~\eqref{eq:main_formulation}.
The perspective formulation is a mixed-integer convex program (MICP), for which specialized solvers with reasonable practical performance exist.
The perspective formulation can also be used to compute a lower bound
on the original objective by relaxing the integral constraints, as in~\cite{moehle2015perspective},
as well to obtain good initializations for any of the procedures described in~\S\ref{sec:methods}.

\paragraph{Perspective.}
Following~\cite[\S 8]{rockafellar1970convex}, we define the perspective (or recession) of the
closed convex function $f$ with $0 \in \dom f$
as\footnote{If $0 \not \in \dom f$, replace $\gamma f_0(x / \gamma)$ with $\gamma f_0(y + x / \gamma)$ for any $y \in \dom f$. See~\cite[Thm.\ 8.3]{rockafellar1970convex} for more details.}
\BEQ
\label{eq:persp}
f^\mathrm{p}(x, t) = \begin{cases}
  t f(x / t) & t > 0,\\
  \lim_{\gamma \downto 0}\,\gamma f_0(x/\gamma) & t = 0,\\
  +\infty & \text{otherwise},
\end{cases}
\EEQ
for $(x, t) \in \reals^n \times \reals_+$.
We will use the fact that the resulting function $f^\mathrm{p}$ is convex~\cite[\S3.2.6]{boyd2004convex}.

\paragraph{Superlinearity assumption.}
If $f$ is superlinear, \ie, if for all $x \in \reals^n \setminus \{0\}$, we have
\BEQ\label{eq:superlinear-assumption}
\lim_{t \to \infty} \frac{f(tx)}{t} = +\infty,
\EEQ
then
\BEQ
\label{eq:superlinear}
f^\mathrm{p}(x, t) = \begin{cases}
  tf(x/t) & t > 0 \\
  0 & t = 0, \; x = 0,\\
  +\infty & \text{otherwise},
 \end{cases}
\EEQ
since the limit in~\eqref{eq:persp} is equal to the limit in~\eqref{eq:superlinear-assumption} unless $x=0$.

There are many convex functions that satisfy this superlinearity property.
Some examples are the sum of squares function and the indicator function of a compact convex set.
Since we will make heavy use of property~\eqref{eq:superlinear} in this section,
we will assume that $f_0$ is superlinear for the remainder of this section.
If $f_0$ is not superlinear, then it can be made superlinear
by adding, \eg, a small positive multiple of the sum of squares function.

\paragraph{Conic representation of the perspective.}
We note that representing the epigraph of the perspective of a function
is often simple if the function has a conic representation~\cite{grant2008graph}.
More specifically, if $f$ has a conic representation
\[
f(x) \le v \iff Ax + bv + c \in \mathcal K,
\]
for some closed convex cone $\mathcal K$, then the perspective of $f$ has a conic representation given by
\[
f^\mathrm{p}(x, t) \le v \iff Ax + bv + tc \in \mathcal K.
\]
This fact allows us to use a conic representation of the perspective
and avoid issues of non-differentiability and division-by-zero that
we might encounter with direct numerical implementations of the perspective~\cite[\S 2]{moehle2015perspective}.

\paragraph{Perspective formulation.}
We define the \emph{perspective formulation} of~\eqref{eq:main_formulation} as the following MICP:
\BEQ
\begin{array}{ll}
\mbox{minimize} & \sum_{i=1}^m f^\mathrm{p}_i(z_i, t_i) + (1 - t_i) \alpha_i +
\frac{1}{m}\left(f^\mathrm{p}_0(z_i, t_i) + f^\mathrm{p}_0(x - z_i, 1-t_i)\right) \\
\mbox{subject to} & t \in \{0, 1\}^m,
\end{array}
\label{eq:micp}
\EEQ
with variables $x, z_i\in\reals^n$ for $i=1,\ldots,m$ and $t\in\reals^m$.
Any MICP solver that can handle the functions $f^\mathrm{p}_i$ for $i=0,\ldots,m$
can be used to solve~\eqref{eq:micp}.

\paragraph{Proof of equivalence.}
To show that~\eqref{eq:micp} is equivalent to the original problem~\eqref{eq:main_formulation},
first take $(x, t, z_i)$ that are feasible for~\eqref{eq:micp}.
Since $t$ is Boolean, for each $i$ we have $t_i=0$ or $t_i=1$.
Since $f^\mathrm{p}_0(z_i,t_i)$ must be finite (as this point is feasible),
then $t_i=0$ implies that $z_i=0$ (due to~\eqref{eq:superlinear}).
Similarly, when $t_i = 1$ we must have $z_i = x$.
Therefore the $i$th term in the sum becomes
\[
t_if_i(x) + (1-t_i)\alpha_i + \frac{1}{m}f_0(x).
\]
Summing over the index $i$ yields that problem~\eqref{eq:micp} is equivalent to
\BEQ
\begin{array}{ll}
\mbox{minimize} & f_0(x) + \sum_{i=1}^m t_i f_i(x) + (1 - t_i)\alpha_i \\
\mbox{subject to} & t \in \{0, 1\}^m.
\end{array}
\label{eq:aux-mip}
\EEQ
Partially minimizing~\eqref{eq:aux-mip} over $t$, we find that $x$ is a feasible
point for~\eqref{eq:main_formulation} with the same objective value.

Now take $x$ feasible for~\eqref{eq:main_formulation}.
Let
\[
t_i = \begin{cases}
1 & f_i(x) \leq \alpha_i \\
0 & \text{otherwise},
\end{cases}
\quad
i=1,\ldots,m,
\]
and $z_i = t_ix$.
Then $(x, t, z_i)$ is feasible for~\eqref{eq:micp} and has the same objective value, and the problems are equivalent.

\paragraph{Lower bound via relaxation.}
Since the perspective formulation is equivalent to the original problem,
relaxing the Boolean constraint in~\eqref{eq:micp}
and solving the resulting convex optimization problem
\BEQ
\begin{array}{ll}
\mbox{minimize} & \sum_{i=1}^m f^\mathrm{p}_i(z_i, t_i) + (1 - t_i) \alpha_i +
\frac{1}{m}\left(f^\mathrm{p}_0(z_i, t_i) + f^\mathrm{p}_0(x - z_i, 1-t_i)\right) \\
\mbox{subject to} & 0 \leq t \leq \ones,
\end{array}
\label{eq:relaxed}
\EEQ
with variables $z_i$, $t$, and $x$,
yields a lower bound on the objective value of~\eqref{eq:main_formulation}.
That is, given any approximate solution of~\eqref{eq:main_formulation}
with objective value $p$, the optimal value $q^\star$ of~\eqref{eq:relaxed}
yields a certificate guaranteeing that the approximate solution
is suboptimal by at most $p-q^\star$.
Additionally, a solution of the relaxed problem can be used as an
initial point for any of the heuristic methods described in~\S\ref{sec:methods}.

\paragraph{Efficiently solving the relaxed problem.}
We note that~\eqref{eq:relaxed} has $m+1$ times as many variables as the original problem, so it is worth considering faster solution methods.
To do so, we can convert the problem to
\emph{consensus form}~\cite[\S7.1]{boyd2011distributed};
\ie, we introduce additional variables $y_i\in\reals^n$ for $i=1,\ldots,m$,
and constrain $y_i=x$, resulting in the equivalent problem
\BEQ
\begin{array}{ll}
\mbox{minimize} & \sum_{i=1}^m f^\mathrm{p}_i(z_i, t_i) + (1 - t_i) \alpha_i +
\frac{1}{m}\left(f^\mathrm{p}_0(z_i, t_i) + f^\mathrm{p}_0(y_i - z_i, 1-t_i)\right) \\
\mbox{subject to} & y_i=x,\quad i=1,\ldots,m,\\
 & 0 \leq t \leq \ones.
\end{array}
\label{eq:persp-relaxation}
\EEQ
Since the objective is separable in $(y_i,z_i,t_i)$ over $i$,
there exist many efficient distributed algorithms for solving this problem, \eg,
the alternating direction method of multipliers (ADMM)
\cite{boyd2011distributed, M2AN_1975__9_2_41_0, gabay1976dual}.

\section{Implementation}
\label{sec:implementation}

Our Python package \verb|sccf| approximately solves generic problems
of the form~\eqref{eq:main_formulation} provided all $f_i$ can be represented as valid \verb|cvxpy| expressions and constraints.
It is available at:
\begin{center}
\texttt{https://www.github.com/cvxgrp/sccf}.
\end{center}
We provide a method \verb|sccf.minimum|, which
can be applied to a \verb|cvxpy| Expression and a scalar to create a
\verb|sccf.MinExpression|.
The user then forms an objective as a sum of \verb|sccf.MinExpression|s,
passes this objective and (possibly) constraints to a \verb|sccf.Problem|
object, and then calls the \verb|solve| method, which
implements algorithm~\ref{alg:twostep}.
We take advantage of the fact that the only parameter changing
between problems is $\lambda$ by caching the canonicalization procedure \cite{agrawal2019differentiable}. 
Here is an example of using \verb|sccf| to solve a clipped least squares problem:
\begin{verbatim}
import cvxpy as cp
import sccf

A, b = get_data(m, n)

x = cp.Variable(n)
objective = 0.0
for i in range(m):
    objective += sccf.minimum(cp.square(A[i]@x-b[i]), 1.0)
objective += 0.01 * cp.sum_squares(x)

prob = sccf.Problem(objective)
prob.solve()
\end{verbatim}

\section{Examples}
\label{sec:examples}

All experiments were conducted on a single core
of an Intel i7-8700K CPU clocked at 3.7 GHz.

\subsection{Clipped regression}
\label{sec:erm-example}
In this example we compare clipped regression (\S\ref{sec:clipped_erm})
with standard linear regression and
Huber regression~\cite{huber1973robust}
(a well known technique for robust regression)
on a one-dimensional dataset with outliers.
We generated data by sampling 20 data points $(x_i,y_i)$ according to
\[
x_i \sim \mathcal N(0,1), \quad y_i = x_i + (0.1)z_i,
\quad z_i\sim\mathcal N(0,1),
\quad i=1,\ldots,20.
\]
We introduced outliers in our data by flipping the sign of $y_i$
for 5 random data points.

The problems all have the form
\BEQ
\begin{array}{ll}
\mbox{minimize} & L(\theta) = \sum_{i=1}^{20}\phi(x_i\theta-y_i) + (0.2)\theta^2,
\end{array}
\label{eq:clipped_regression}
\EEQ
where $\phi:\reals\to\reals$ is a penalty function.
In clipped regression, $\phi(z) = \min\{z^2, 0.5\}$.
In linear regression, $\phi(z) = z^2$.
In Huber regression,
\[
\phi(z) = \begin{cases}
z^2 & |z| \leq 0.5 \\
0.5(2|z| - 0.5) & \text{otherwise}.
\end{cases}
\]

\begin{figure}
    \centering
    \includegraphics[width=.6\linewidth]{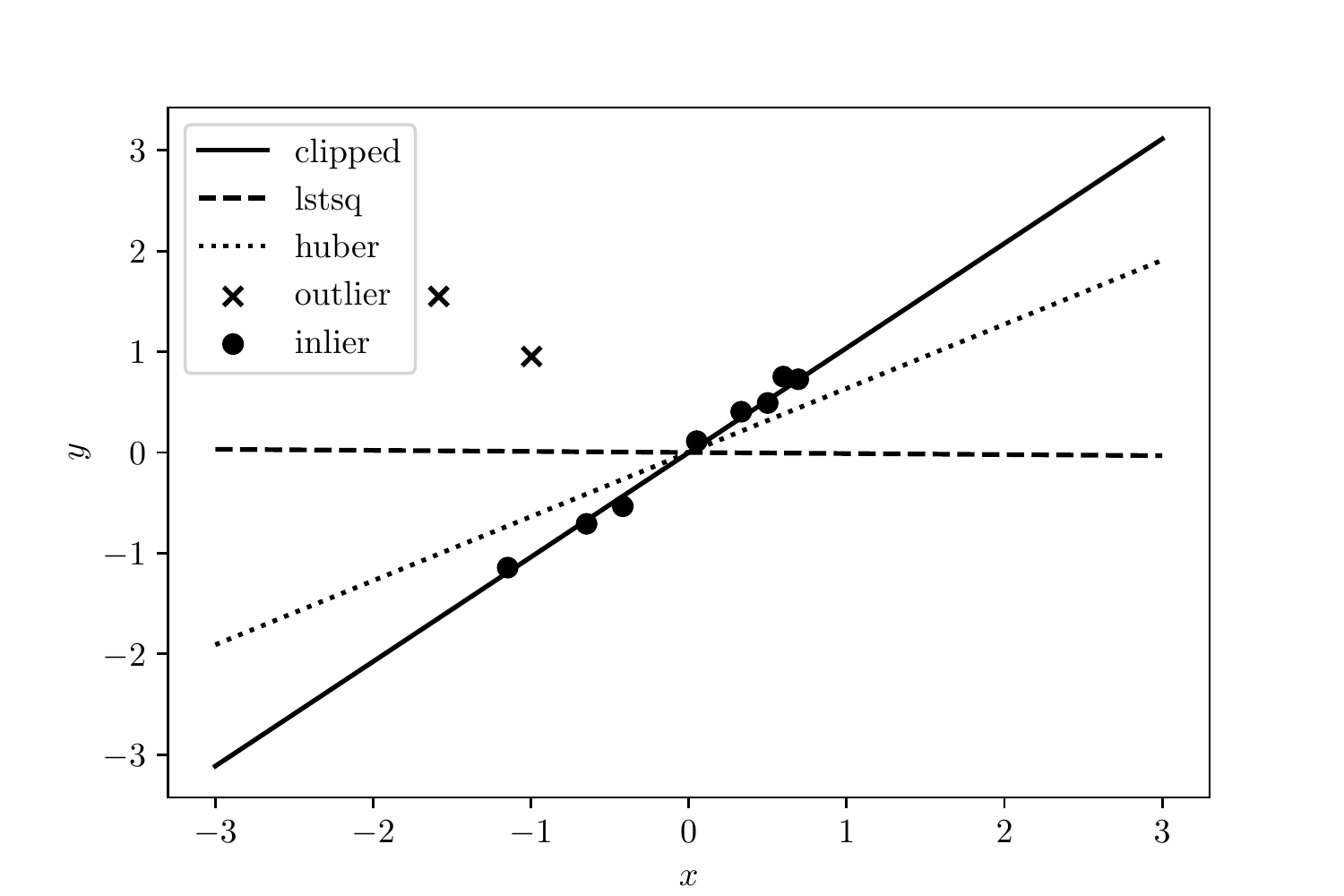}
    \caption{Clipped regression, linear regression, and Huber regression on a one-dimensional dataset with outliers. The outliers affect the linear regression and Huber regression models, while the clipped regression model appears to be minimally affected.}
    \label{fig:clipped-regression}
\end{figure}

Let $\theta^\mathrm{clip}$ be the clipped regression model;
we deem points where $(x_i\theta^\mathrm{clip}-y_i)^2\geq 0.5$ as outliers
and the remaining points as inliers.
In figure~\ref{fig:clipped-regression} we visualize
the data points and the resulting models along with the outliers/inliers
identified by the clipped regression model.
In this figure, the clipped regression model clearly outperforms the linear and Huber
regression models since it is able to fully ignore the outliers.
Algorithm~\ref{alg:twostep} terminated in 0.13 seconds and took 8 iterations on this instance. 

\paragraph{Lower bound.}
The relaxed version of the perspective formulation~\eqref{eq:relaxed}
can be used to efficiently find a lower bound on the objective value for the clipped version
of~\eqref{eq:clipped_regression}. The objective value of~\eqref{eq:clipped_regression} for clipped regression was 1.147,
while the lower bound we calculated was 0.533,
meaning our approximate solution is suboptimal by at most 0.614.

\begin{figure}
    \centering
    \includegraphics[width=.7\linewidth]{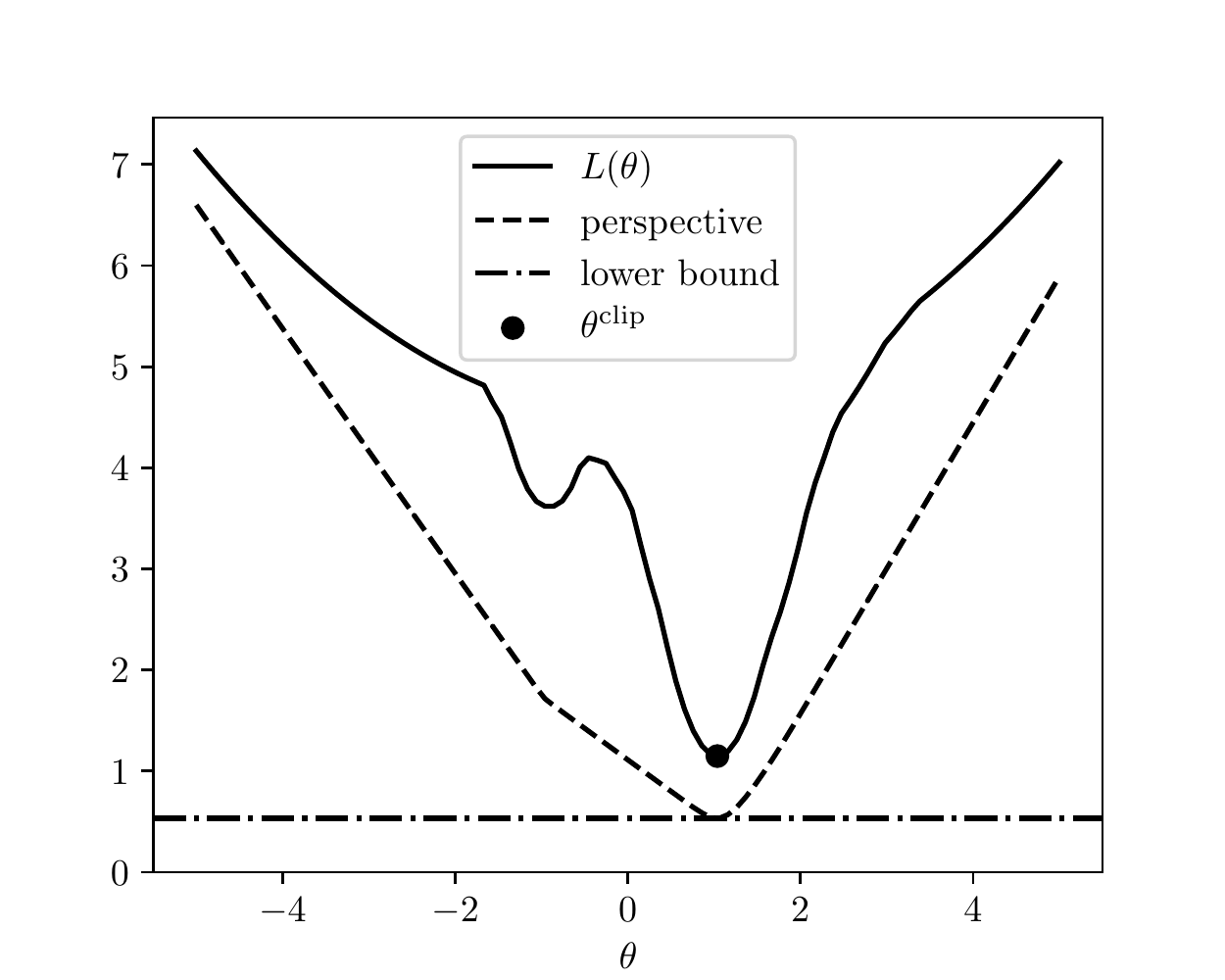}
    \caption{The clipped regression loss and its perspective relaxation.}
    \label{fig:perspective}
\end{figure}

In figure~\ref{fig:perspective} we plot the clipped objective~\eqref{eq:clipped_regression}
for various values of $\theta$;
note that the function is highly nonconvex and that $\theta^\mathrm{clip}$
is the (global) solution.
We also plot the objective of the perspective relaxation as a function of $\theta$,
found by partially minimizing~\eqref{eq:relaxed} over $z_i$ and $t$; note that
the function is convex and a surprisingly good approximation of the true convex envelope.
We note that the minimum of the perspective relaxation
and the true minimum are surprisingly close, leading us to believe that
the solution to the perspective relaxation could be a good
initialization for heuristic methods.

\subsection{Clipped logistic regression}
In this example we apply clipped logistic regression (\S\ref{sec:clipped_erm})
to a dataset with outliers.
We generated data by sampling 1000 data points $(x_i,y_i)$
from a mixture of two Gaussian distributions in $\reals^5$.
We randomly partitioned the data into 100 training data points
and 900 test data points and introduced outliers by flipping the sign of $y_i$
for 20 random training data points.

We (approximately) solved the \emph{clipped logistic regression} problem
\[
\begin{array}{ll}
\mbox{minimize} & \frac{1}{1000}
\sum_{i=1}^{1000}\min\{\log(1 + e^{-y_i (x_i^T\theta + b)}),\alpha\} + (0.1)\|\theta\|_2^2,
\end{array}
\]
with variables $\theta$ and $b$, for various values of $\alpha\in[10^{-1},10^1]$.
We also solved the problem for $\alpha=+\infty$, \ie, the \emph{standard logistic
regression problem}.
Over the $\alpha$ values we tried, on average,
algorithm~\ref{alg:twostep} took 6.37 seconds
and terminated in 9.64 iterations.

\begin{figure}
    \centering
    \includegraphics[width=.7\linewidth]{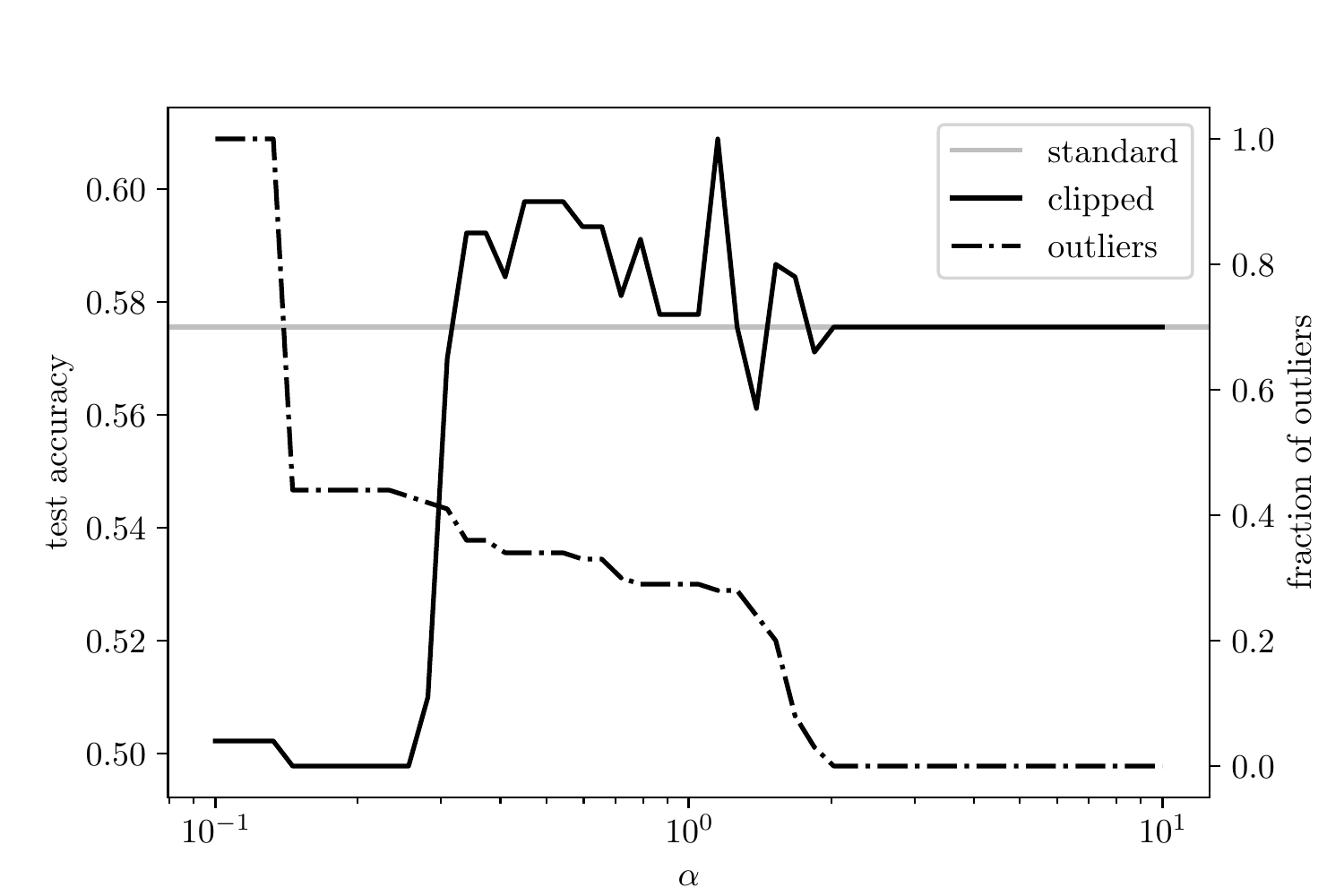}
    \vspace{-1em}
    \caption{Test accuracy of clipped logistic regression (solid),
    test accuracy of standard logistic regression (gray),
    and fraction of outliers (dotted dashed) for varying clip values $\alpha$.
    Note that the fraction of detected outliers goes down as $\alpha$ goes up.
    Between roughly $\alpha=10^{-.5}$ and $\alpha=10^{0.05}$, the test
    accuracy of clipped
    logistic regression is higher than standard logistic regression.
    Clipped logistic regression converges to standard logistic regression as
    $\alpha \to \infty$.
    }
    \label{fig:logreg}
\end{figure}

\begin{figure}
    \centering
    \includegraphics[width=.6\linewidth]{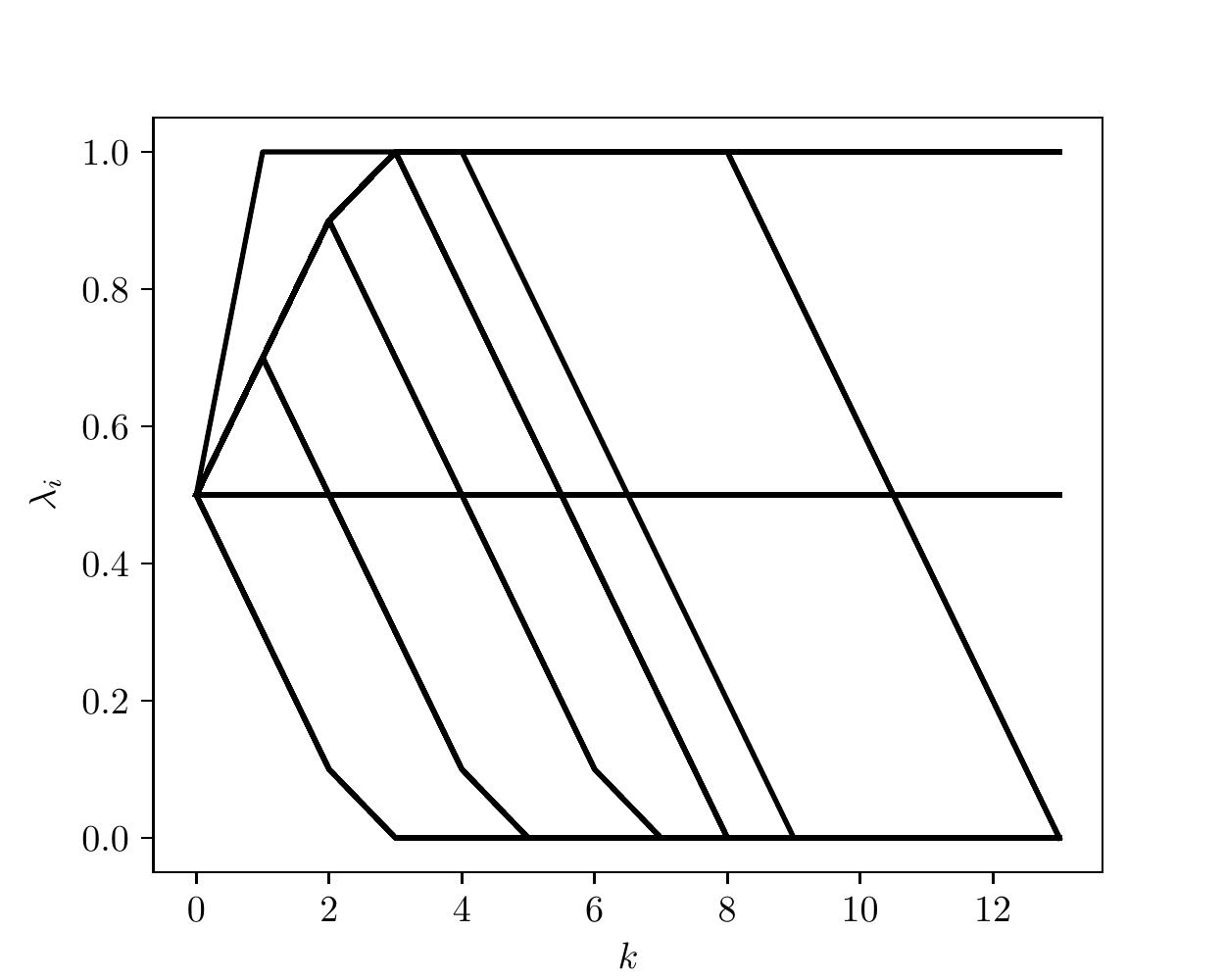}
    \vspace{-1em}
    \caption{A plot of $\lambda$ throughout the course of algorithm~\ref{alg:twostep}
    for the clipped logistic regression example. Note that at some of the iterations
    (\eg, $k=1$, $2$, or $3$),
    the gradient of the loss with respect to a certain $\lambda_i$ changes sign,
    causing $\lambda_i$ to be updated in the opposite direction.}
    \label{fig:logreg_lambda}
\end{figure}

Figure~\ref{fig:logreg} displays the test loss and fraction of outliers
over the range of values of $\alpha$ we approximately minimized.
Figure~\ref{fig:logreg_lambda} shows the trajectory of the entries of $\lambda$
during each step of the execution of algorithm~\ref{alg:twostep} for the $\alpha$ with the highest test accuracy, while figure~\ref{fig:logreg_lambda} plots the histogram of the logistic loss for each of the available data points for this same $\alpha$.

\begin{figure}
    \begin{subfigure}[b]{0.4\textwidth}
        \centering
        \includegraphics[]{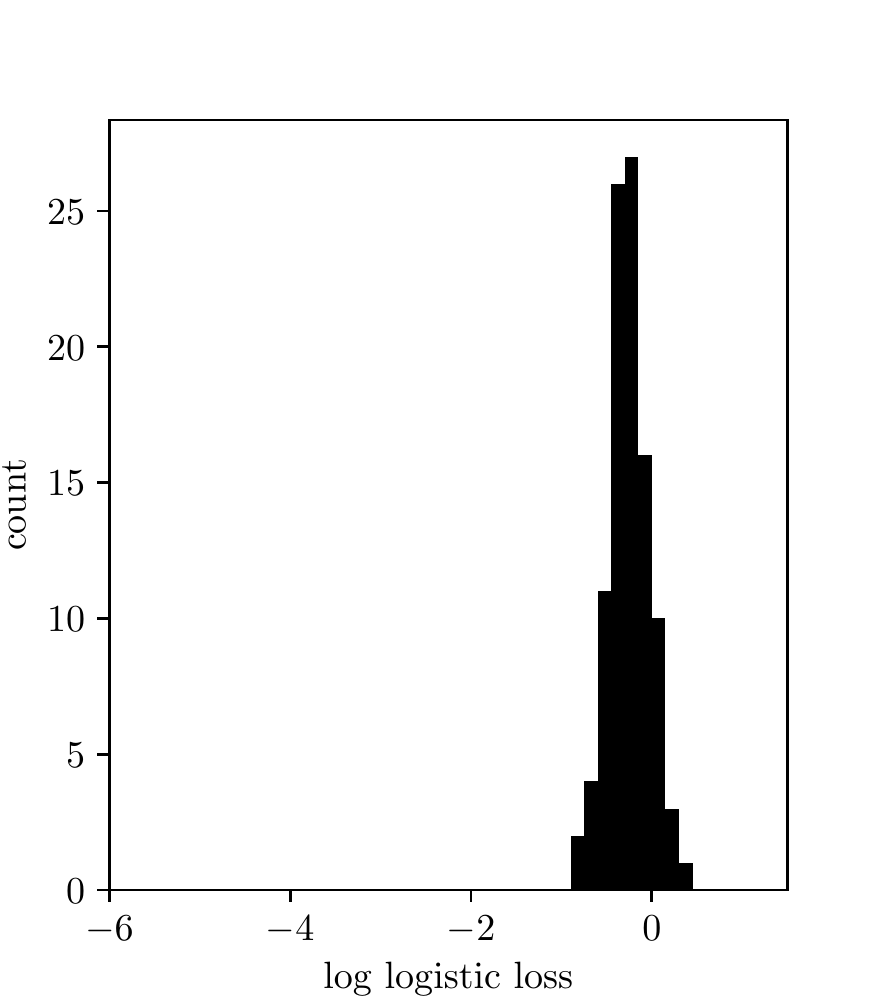}
    \end{subfigure}
    \begin{subfigure}[b]{0.7\textwidth}
        \centering
        \includegraphics[]{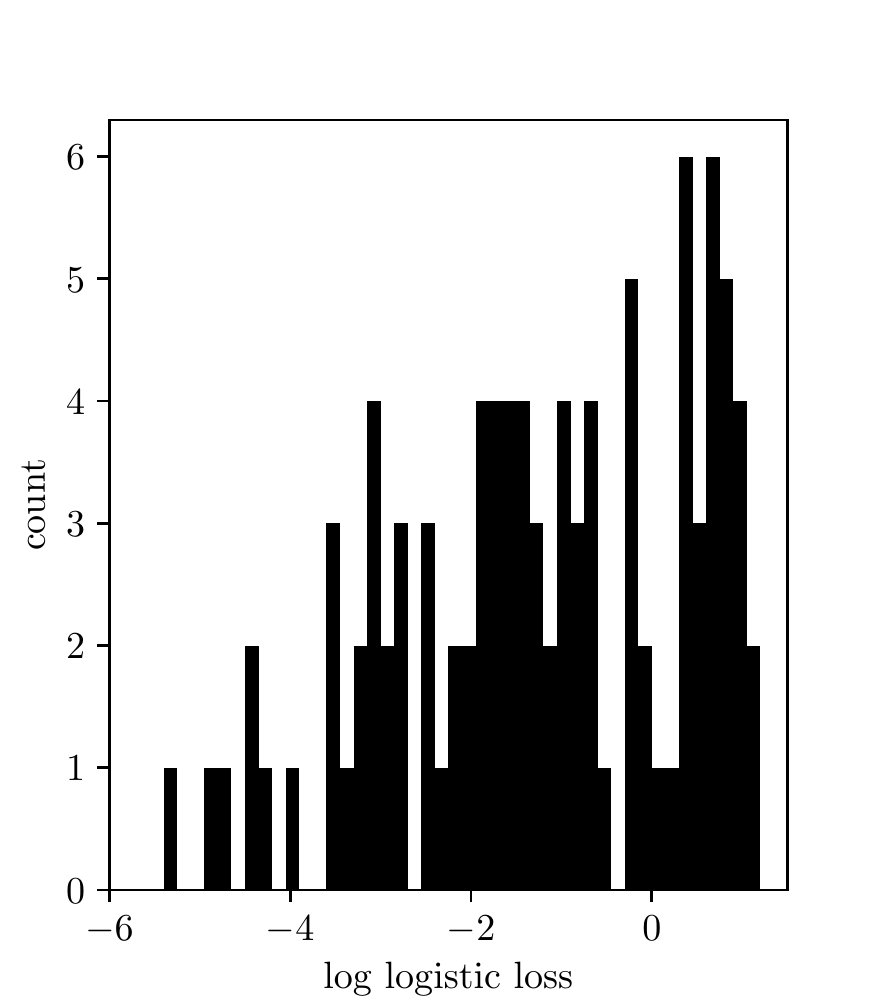}
    \end{subfigure}
    \vspace{-1.5em}
    \caption{Left: histogram of log logistic loss for each
    data point in standard
    logistic regression; right: histogram of log logistic loss
    for each data point in clipped logistic regression. Note that standard logistic
    regression attempts to make the loss small for all data points,
    while its clipped counterpart allows the loss to be high for some of the data points.}
    \label{fig:density}
\end{figure}

\subsection{Lane changing}
\label{sec:control-example}

In this example, we consider a control problem
where a vehicle traveling down a road at a fixed speed
must avoid obstacles, stay in one of two lanes,
and provide a comfortable ride.
We let $x_t\in\reals$ denote the lateral position of the vehicle
at time $t=0,\ldots,T$ ($T$ is the time horizon).

The obstacle avoidance constraints are given as vectors
$x^\mathrm{min}, x^\mathrm{max}\in\reals^T$ that
represent lower and upper bounds on $x_t$ at time $t$.

We can split the objective into the sum of two functions described below.
\begin{itemize}
    \item \emph{Lane cost.} Suppose the two lanes are centered at $x=-1$ and $x=1$.
    The lane cost is given by 
    \[
    g^\mathrm{lane}(x) = \sum_{t=0}^T \min\{(x_t-1)^2,1\} + \min\{(x_t+1)^2,1\}.
    \]
    The lane cost incentivizes the vehicle to be in the center of
    one of the two lanes.
    The lane cost is evidently a sum of clipped convex functions.
    \item \emph{Comfort cost.}
    The comfort cost is given by
    \[
    g^\mathrm{comfort}(x) = \rho_1 \|Dx\|_2^2 + \rho_2 \|D^2x\|_2^2 + \rho_3 \|D^3x\|_2^2,
    \]
    where $D$ is the difference operator and $\rho_1,\rho_2,\rho_3 > 0$ are weights
    to be chosen.
    The comfort cost is a weighted sum of the squared lateral velocity, acceleration, and jerk.
\end{itemize}
To find the optimal lateral trajectory we solve the problem
\begin{equation}
\begin{array}{ll}
\mbox{minimize} & g^\mathrm{lane}(x) + g^\mathrm{comfort}(x) \\
\mbox{subject to} & x_0=x^\mathrm{start}, \quad x_T=x^\mathrm{end}, \\
& x^\mathrm{min} \leq x \leq x^\mathrm{max},
\end{array}
\label{eq:clipped-control-example}
\end{equation}
where $x^\mathrm{start},x^\mathrm{end}\in\reals$ are given starting and ending
points of the trajectory.

\begin{figure}
    \centering
    \includegraphics[]{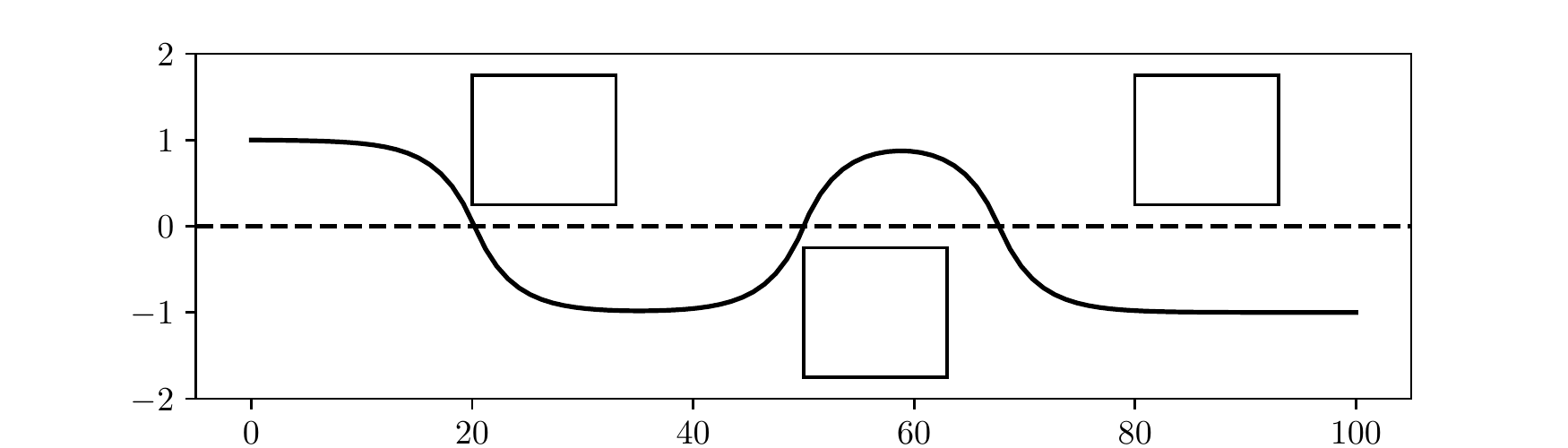}
    \caption{Trajectory of a vehicle looking to avoid obstacles
    (represented by boxes) while optimizing for comfort and lane position.}
    \label{fig:clipped-control-complex}
\end{figure}

\paragraph{Numerical example.}
We use $T=100$, $\rho_1 = 10$, $\rho_2 = 1$, $\rho_3 = .1$,
$x^\mathrm{start}=1$, and $x^\mathrm{end}=-1$.
In figure~\ref{fig:clipped-control-complex}
we show the trajectory resulting from an approximate solution
to~\eqref{eq:clipped-control-example} with three obstacles.
For this example, algorithm~\ref{alg:twostep} terminated in 1.2 seconds
and took 4 iterations.
We are able to find a comfortable trajectory that avoid the obstacles
and spends as little time as possible in between the lanes.

\paragraph{Lower bound.}
Using the relaxed version of the perspective formulation~\eqref{eq:relaxed}, we can compute a lower bound on the objective value of the clipped control problem~\eqref{eq:clipped-control-example}. We found a lower bound value of around 103.55, while the approximate solution we found had an objective value of 119.07, indicating that our approximate solution is no more
than 15\% suboptimal.

\section*{Acknowledgments}
S. Barratt is supported by the National Science Foundation Graduate Research Fellowship
under Grant No. DGE-1656518.

\bibliography{refs}

\begin{thebibliography}{10}

\bibitem{agrawal2019differentiable}
A.~Agrawal, B.~Amos, S.~Barratt, S.~Boyd, S.~Diamond, and Z.~Kolter.
\newblock Differentiable convex optimization layers.
\newblock In {\em Advances in Neural Information Processing Systems}, 2019.

\bibitem{boyd2011distributed}
S.~Boyd, N.~Parikh, E.~Chu, B.~Peleato, and J.~Eckstein.
\newblock Distributed optimization and statistical learning via the alternating
  direction method of multipliers.
\newblock {\em Foundations and Trends{\textregistered} in Machine Learning},
  3(1):1--122, 2011.

\bibitem{boyd2004convex}
S.~Boyd and L.~Vandenberghe.
\newblock {\em Convex optimization}.
\newblock Cambridge University Press, 2004.

\bibitem{gabay1976dual}
D.~Gabay and B.~Mercier.
\newblock A dual algorithm for the solution of nonlinear variational problems
  via finite element approximation.
\newblock {\em Computers \& Mathematics with Applications}, 2(1):17--40, 1976.

\bibitem{M2AN_1975__9_2_41_0}
R.~Glowinski and A.~Marroco.
\newblock Sur l'approximation, par \'el\'ements finis d'ordre un, et la
  r\'esolution, par p\'enalisation-dualit\'e d'une classe de probl\`emes de
  dirichlet non lin\'eaires.
\newblock {\em ESAIM: Mathematical Modelling and Numerical Analysis -
  Mod\'elisation Math\'ematique et Analyse Num\'erique}, 9(R2):41--76, 1975.

\bibitem{grant2008graph}
M.~Grant and S.~Boyd.
\newblock Graph implementations for nonsmooth convex programs.
\newblock In {\em Recent Advances in Learning and Control}, pages 95--110.
  Springer, 2008.

\bibitem{huber1973robust}
P.~Huber.
\newblock Robust regression: asymptotics, conjectures and monte carlo.
\newblock {\em The Annals of Statistics}, 1(5):799--821, 1973.

\bibitem{huber2009robust}
P.~Huber and E.~Ronchetti.
\newblock {\em Robust Statistics}.
\newblock John Wiley \& Sons, 2009.

\bibitem{lan2016robust}
G.~Lan, C.~Hou, and D.~Yi.
\newblock Robust feature selection via simultaneous capped ℓ2-norm and ℓ2,
  1-norm minimization.
\newblock In {\em IEEE Intl. Conf. on Big Data Analysis (ICBDA)}, pages 1--5.
  IEEE, 2016.

\bibitem{lipp2016variations}
T.~Lipp and S.~Boyd.
\newblock Variations and extension of the convex--concave procedure.
\newblock {\em Optimization and Engineering}, 17(2):263--287, 2016.

\bibitem{minimizing2019liu}
T.~Liu and H.~Jiang.
\newblock Minimizing sum of truncated convex functions and its applications.
\newblock {\em Journal of Computational and Graphical Statistics}, 28(1):1--10,
  2019.

\bibitem{moehle2015perspective}
N.~Moehle and S.~Boyd.
\newblock A perspective--based convex relaxation for switched-affine optimal
  control.
\newblock {\em Systems \& Control Letters}, 86:34--40, 2015.

\bibitem{ong2013learning}
C.~Ong and L.~An.
\newblock Learning sparse classifiers with difference of convex functions
  algorithms.
\newblock {\em Optimization Methods and Software}, 28(4):830--854, 2013.

\bibitem{portilla2015efficient}
J.~Portilla, A.~Tristan-Vega, and I.~Selesnick.
\newblock Efficient and robust image restoration using multiple-feature
  l2-relaxed sparse analysis priors.
\newblock {\em IEEE Transactions on Image Processing}, 24(12):5046--5059, 2015.

\bibitem{rockafellar1970convex}
T.~Rockafellar.
\newblock {\em Convex analysis}.
\newblock Princeton University Press, 1970.

\bibitem{safari2014insensitive}
A.~Safari.
\newblock An e--{E}--insensitive support vector regression machine.
\newblock {\em Computational Statistics}, 29(6):1447--1468, 2014.

\bibitem{she2011outlier}
Y.~She and A.~Owen.
\newblock Outlier detection using nonconvex penalized regression.
\newblock {\em Journal of the American Statistical Association},
  106(494):626--639, 2011.

\bibitem{sun2013robust}
Q.~Sun, S.~Xiang, and J.~Ye.
\newblock Robust principal component analysis via capped norms.
\newblock In {\em Proc. Intl. Conf. on Knowledge Discovery and Data Mining},
  pages 311--319. ACM, 2013.

\bibitem{suzumura2014outlier}
S.~Suzumura, K.~Ogawa, M.~Sugiyama, and I.~Takeuchi.
\newblock Outlier path: A homotopy algorithm for robust {SVM}.
\newblock In {\em Intl. Conf. on Machine Learning}, pages 1098--1106, 2014.

\bibitem{tao1997convex}
P.~Tao and L.~An.
\newblock Convex analysis approach to {DC} programming: Theory, algorithms and
  applications.
\newblock {\em Acta Mathematica Vietnamica}, 22(1):289--355, 1997.

\bibitem{torr1998robust}
P.~Torr and A.~Zisserman.
\newblock Robust computation and parametrization of multiple view relations.
\newblock In {\em Intl. Conf. on Computer Vision}, pages 727--732. IEEE, 1998.

\bibitem{xu2016robust}
G.~Xu, B.-G. Hu, and J.~Principe.
\newblock Robust {C}-loss kernel classifiers.
\newblock {\em IEEE Transactions on Neural Networks and Learning Systems},
  29(3):510--522, 2016.

\bibitem{yang2010relaxed}
Y.-l. Yu, M.~Yang, L.~Xu, M.~White, and D.~Schuurmans.
\newblock Relaxed clipping: A global training method for robust regression and
  classification.
\newblock In {\em Advances in Neural Information Processing Systems}, pages
  2532--2540, 2010.

\bibitem{yuille2003concave}
A.~Yuille and A.~Rangarajan.
\newblock The concave--convex procedure.
\newblock {\em Neural Computation}, 15(4):915--936, 2003.

\bibitem{zhang2009multi}
T.~Zhang.
\newblock Multi-stage convex relaxation for learning with sparse
  regularization.
\newblock In {\em Advances in Neural Information Processing Systems}, pages
  1929--1936, 2009.

\bibitem{zhang2010analysis}
T.~Zhang.
\newblock Analysis of multi-stage convex relaxation for sparse regularization.
\newblock {\em Journal of Machine Learning Research}, 11(Mar):1081--1107, 2010.

\end{thebibliography}

\appendix

\section{Difference of convex formulation}
\label{sec:convex_concave}
In this section we make the observation that~\eqref{eq:main_formulation} can be expressed as a difference of convex (DC) programming
problem.

Let $h_i(x) = \max(f_i(x) - \alpha_i, 0)$.
This (convex) function measures how far $f_i(x)$ is above $\alpha_i$.
We can express the $i$th term in the sum as
\[
\min\{f_i(x), \alpha_i\} = f_i(x) - h_i(x),
\]
since when $f_i(x) \leq \alpha_i$, we have $h_i(x)=0$,
and when $f_i(x) > \alpha$, we have $h_i(x)=f_i(x)-\alpha_i$.
Since $f_i$ and $h_i$ are convex,~\eqref{eq:main_formulation}
can be expressed as the DC programming problem
\BEQ
\begin{array}{ll}
\mbox{minimize} & f_0(x) + \sum_{i=1}^m f_i(x) - \sum_{i=1}^m h_i(x),
\end{array}
\label{eq:dc-prob}
\EEQ
with variable $x$.
We can apply then well-known algorithms like the convex-concave procedure~\cite{tao1997convex,yuille2003concave}
to (approximately) solve~\eqref{eq:dc-prob}.

\section{Minimal convex extension}

If we replace each $f_i$ with any function $\tilde f_i$ such that
$\tilde f_i(x) = f_i(x)$ when $f_i(x) \leq \alpha_i$,
we get an equivalent problem.
One such $\tilde f_i$ is the \emph{minimal convex extension} of $f_i$,
which is given by
\[
\tilde f_i(x) \coloneqq \sup \{f_i(z) + g^T (x-z) \mid g \in \partial f_i(z), f_i(z) \leq \alpha_i, z\in\reals^n\}.
\]
In general, the minimal convex extension of a function is often hard to compute, but
it can be represented analytically in some (important) special cases.
For example, if $f_i(x)=(a^T x - b)^2$, the minimal convex extension is
the Huber penalty function, or
\[
\tilde f_i(x) = \begin{cases}
(a^Tx - b)^2 & |a^Tx-b| \leq \alpha_i \\
\alpha_i (2|a^Tx-b| - \alpha_i) & \text{otherwise}.
\end{cases}
\]
Using the minimal convex extension leads to an equivalent problem, but,
depending on the algorithm, replacing $f_i$ with $\tilde f_i$ can lead to better numerical performance.

\end{document}